\documentclass [11pt]{amsart}
\usepackage {amsmath, amssymb, a4wide, epsfig, enumerate}
\usepackage[latin1]{inputenc}
\date{March 2005}

\author{Romain Dujardin}
\title{Laminar currents and  birational dynamics}

\newcommand{\cc}{\mathbb{C}}
\newcommand{\bb}{\mathbb{B}}

\newcommand{\dd}{\mathbb{D}}
\newcommand{\zz}{\mathbb{Z}}
\newcommand{\nn}{\mathbb{N}}
\newcommand{\e}{\varepsilon}
\newcommand{\cv}{\rightarrow}

\newcommand{\fr}{\partial}
\newcommand{\om}{\Omega}
\newcommand{\set}[1]{\left\{#1\right\}}
\newcommand{\norm}[1]{\left\Vert#1\right\Vert}
\newcommand{\abs}[1]{\left\vert#1\right\vert}

\newcommand{\cd}{\cc^2}
\newcommand{\pd}{{\mathbb{P}^2}}
\newcommand{\pu}{{\mathbb{P}^1}}
\newcommand{\rest}[1]{ \arrowvert_{#1}}
\newcommand{\m}{{\bf M}}

\newcommand{\unsur}[1]{\frac{1}{#1}}
\newcommand{\el}{\mathcal{L}}

\newcommand{\qq}{\mathcal{Q}}
\newcommand{\geom}{\hspace{.3em}\dot{\wedge}\hspace{.3em}}

\DeclareMathOperator{\area}{Area}

\DeclareMathOperator{\supp}{Supp}
\DeclareMathOperator{\vol}{Vol}

\newtheorem{prop} {Proposition} [section]
\newtheorem{thm}[prop] {Theorem} 
\newtheorem{defi}[prop] {Definition}
\newtheorem{lem}[prop] {Lemma}

\newtheorem{theo}{Theorem} 
\newtheorem{exam}[prop]{Example}

\newtheorem{rmk}[prop]{Remark}

\newenvironment{pf}{\noindent {\bf Proof:} }{\hfill $\square$\\}



\begin{document}

\begin{abstract} 
We study the dynamics of a bimeromorphic map $X\rightarrow X$, where $X$ is a 
compact complex K{\"a}hler surface. Under  a natural geometric hypothesis, 
we construct
an invariant probability measure, which is mixing, hyperbolic and of maximal entropy.
The proof relies heavily on the theory of laminar currents and 
is new even in the case 
of polynomial automorphisms of $\mathbb{C}^2$. This extends
recent results by E. Bedford and J. Diller.
\end{abstract}
\maketitle

\section{Introduction}

Let $X$ be a compact complex surface and $f$ be a bimeromorphic self
map on $X$. We moreover assume $X$ is K{\"a}hler. 
 We are interested in the study of $(X,f)$ as a dynamical
system. These mappings  generalize
polynomial automorphisms of  $\cd$ (viewed as birational on $\pd$),
 whose dynamics have  turned out to be very rich. The general setting 
 raises interesting problems, both in dynamics and 
in intersection theory of positive
closed currents.  \\

It is now classical to introduce the {\em dynamical degree} 
$\lambda$, which is the asymptotic growth rate of the volumes of 
iterated submanifolds. This number is conjecturally related to
the  topological entropy of $f$ by the equation $h_{top}(f)=\log
\lambda$ (see e.g. V. Guedj \cite{gu} for a general account). 
 It is to be mentioned  that in our context, 
this equality has been subject to intensive numerical study (N. Abarenkova
et al. [Ab1-3]) motivated by questions in statistical physics.\\

An important contribution to the study of the dynamical system $(X,f)$ 
was made by J. Diller and C. Favre \cite{df}. They proved that
the mappings with interesting dynamics are those with $\lambda>1$. Under
this hypothesis, they constructed  positive closed currents $T^\pm$
such that $(f^{\pm 1})^*T^\pm=\lambda T^\pm$. A classical additional
observation is that if $f$ is not birationally conjugate to an
automorphism, then $X$ is a rational surface. 

For the purpose of extending the known results for polynomial
automorphisms, a natural approach is to give a meaning to the
intersection measure $\mu=T^+\wedge T^-$, which should have remarkable
dynamical properties (see e.g. \cite{fg, ca, di2}).
 In the most general context, 
this method, combining pluripotential theory for
the definition of $\mu$, and Pesin's theory for its fine dynamical
study, brings up
several difficulties. The reason is
the presence of indeterminacy points with possibly complicated
dynamics.  
A recent breakthrough is the paper by E. Bedford and J. Diller \cite{bd} in
which they construct the wedge product measure $\mu$ and prove it to
be mixing and hyperbolic (non zero Lyapounov exponents) under the
hypothesis 
\begin{equation}\label{eq_bd1}
\sum_{n\geq 0} \unsur{\lambda^n} \abs{\log {\rm dist}(f^n(I(f^{-1})),
I(f))}<\infty.   
\end{equation}

Our approach 
differs crucially from the previous ones by the systematic use of
the {\em laminar structure} of the currents $T^\pm$. This 
concept dates back to  D. Ruelle and D. Sullivan \cite{rs} 
and was developed by E. Bedford, M. Lyubich and J. Smillie in
their seminal paper \cite{bls}. 

Using the laminar structure allows us to
define the measure $\mu$ without appealing to pluripotential theory,
as the {\em geometric intersection} $\mu=T^+\geom T^-$ of the disks
subordinate to $T^+$ and $T^-$. Next, we derive the dynamical
properties of $\mu$ without use of Pesin's theory, by using an argument
in the style of  M. Lyubich \cite{ly} and 
J.Y. Briend and J. Duval \cite{brd} along the laminar currents. 
The method we use is new even for complex H{\'e}non mappings, and
 provides a new approach for the geometric analysis of
the maximal entropy measure in \cite{bls}, sections 4, 8, and 9. 
Since the Briend-Duval argument also works in higher dimensions 
this approach might open the way to the finer study of $\mu$ in higher
dimension (cf. \cite{dis}).

We now state a precise result. The meaning of the ``algebraically
stable'' assumption in the theorem will be made precise in the next
section. This does not restrict the scope of the theorem,
 since any birational map is
birationally conjugate to an algebraically stable map. 

\begin{theo}\label{theo_entropy}
Let $f$ be an algebraically stable birational map of a rational
surface $X$ with dynamical degree $\lambda>1$.
Assume that the currents $T^+$ and $T^-$ have nontrivial geometric
intersection $\mu= T^+\geom T^-$. Then
\begin{enumerate}[\em i.]
\item $\mu$ is an invariant measure which is mixing.
\item  For $\mu$-almost every $p$, there exist unit tangent vectors
$e^u(p)$ and $e^s(p)$ at $p$, there exists $\nn'\subset\nn$ of density $1$,
such that 
$$
\liminf_{\nn'\ni n\cv\infty} \unsur{n} \log\abs{df^n(e^u(p))}\geq 
\frac{\log\lambda}{2} \text{ and } 
\limsup_{\nn'\ni n\cv\infty} \unsur{n} \log\abs{df^n(e^s(p))}\leq 
-\frac{\log\lambda}{2}.
$$
These bounds are sharp.
\item $\mu$ has entropy $h_\mu(f)=\log\lambda$. 
In particular the topological entropy
 $h_{top}(f)$ is $\log\lambda$.
\item $\mu$ has product structure
with respect to local stable and unstable manifolds. In particular
$(f,\mu)$ has the Bernoulli property.
\end{enumerate}
\end{theo}

In the case of a projective 
non rational surface $X$, $f$ is conjugate to an automorphism on a
torus or $K3$ surface and the result is already known and due to S. Cantat \cite{ca}.
Nevertheless our proof can be adapted so as to apply in this setting
as well. 

From this theorem, it is natural to look for criteria ensuring that 
 $T^+\geom T^->0$. Our second result is the following.

\begin{theo}\label{theo_bd}
Under the assumptions of Theorem \ref{theo_entropy}, assume further
that the Bedford-Diller condition  (\ref{eq_bd1}) holds. Then
 $T^+\geom T^-=T^+\wedge T^->0$, hence Theorem \ref{theo_entropy} applies. Moreover  
$\mu$ describes the asymptotic distribution of saddle orbits, and most
saddle points lie inside $\supp \mu$. 
\end{theo}

Here ``most'' means the following: $f$ has approximately $\lambda^n$
 periodic points of period $n$, and asymptotically (as $n\cv\infty$)
 the number of saddle points inside $\supp \mu$ is equivalent to $\lambda^n$.
 The proof uses classical intersection theory of positive closed currents. 
It would be interesting in view
of getting rid of hypothesis  (\ref{eq_bd1}) to find a completely geometric
argument ensuring that $T^+\geom T^->0$.\\

The structure of the paper is as follows. In \S \ref{sec_structure} we
recall some facts on birational dynamics, mainly from \cite{df}, and 
some results on laminar currents from \cite{structure}
that are crucial in the following. 
In \S \ref{sec_equidist} we prove an
equidistribution property for preimages of points along the unstable
current $T^-$, which is the analogue of the Lyubich-Briend-Duval
lemma \cite{ly,brd}  in our setting. Theorem \ref{theo_entropy} is
proved in in  \S \ref{sec_measure}  and  theorem \ref{theo_bd} 
in \S \ref{sec_bd}.
Another approach to these results, allowing the use of Pesin's theory,
is outlined in the Appendix.  \\

The author would like to thank Eric Bedford for focusing 
his interest on this problem, as well as Jeffrey Diller, Vincent
Guedj, and the anonymous referees
for many constructive  comments.

 %%%%%%%%%%%%%%%%%%%%%%%%%%%%%%%%%%%%%%
 %%%%%%%%%%%%%%%%%%%%%%%%%%%%%%%%%%%%%%
 
\section{Laminar structure}\label{sec_structure}
 
 We first briefly introduce the dynamical setting we consider throughout the paper.
   For more details, the reader is referred to \cite{df,bd}. 

Let $f:X\cv X$ be a bimeromorphic map of a compact K{\"a}hler
surface, with K{\"a}hler form $\omega$. 
We denote by $I(f)$ the indeterminacy set, which is a finite
number of points, and by $C(f)= f^{-1}(I(f^{-1}))$ the critical set. 
We will often use the fact that $f(C(f))=I(f^{-1})$ and $f(I(f))=C(f^{-1})$.

We now review some results of \cite{df}.
First, up to a
bimeromorphic change of surface, we may assume that $f$ is
``algebraically stable'', which means  that 
$$\forall n,m\geq 0, ~ I(f^n)\cap I(f^{-m})=\emptyset. $$ In this case the
dynamical degree
$\lambda$ is the spectral radius of the action of $f^*$ on
$H^{1,1}(X)$. 

From now on we will assume $\lambda>1$.
If $f$ is not birationally conjugate to
an automorphism then $X$ is a rational surface 
(thus our bimeromorphic maps are rather {\em birational}).
The  case of automorphisms of projective non rational surfaces
is treated in \cite{ca}, so  we will assume 
 $X$ is a rational surface (see however Remark \ref{rmk_k3} below).
There exist nef cohomology classes $\theta^+$ and $\theta^-$ in
$H^{1,1}(X)$ such that 
$$\unsur{\lambda}f^*\theta^+ = \theta^+ \text{ and }  
\unsur{\lambda}(f^{-1})^*\theta^- = \theta^-.$$ Moreover there exist
positive closed currents $T^{+/-}$, respectively cohomologous to 
$\theta^{+/-}$, so that for any smooth closed (1,1) form $\alpha$ on $X$, 
$$\lambda^{-k} (f^k)^* \alpha \underset{k\cv\infty}{\cv}
\frac{(\set{\alpha},{\theta^-})}{(\theta^+, \theta^-)} T^+,$$ with a
similar formula for $T^-$, where $\set{\alpha}$ is the cohomology
class of $\alpha$ and $(\cdot,\cdot)$ is the intersection pairing in
cohomology.
 
 Here are some known properties of the currents $T^{+/-}$:
 \begin{itemize}
 \item[-] $f^*(T^+)=\lambda T^+$ and $(f^{-1})^*T^-=\lambda T^-$;
 \item[-] $T^{+/-}$ are extremal in the cone of positive closed currents;
 \item[-] if the Lelong number $\nu(p, T^+)>0$ then $p\in I(f^n)$ for some $n$
 (similarly for $T^-$). In  particular $T^{+/-}$ give no mass to analytic curves. We call 
 such currents {\em diffuse}.
 \end{itemize}
 
\bigskip

The methods in the present article rely very much on some results on 
 the laminar structure of the invariant currents, that we obtained in 
a series of papers [Du1--3]. The motivation was precisely to extend
 the results of \cite{bls} to the widest possible context. We will 
 spend some time to recall some background  on the topic. 

The starting point is the following definition. The definition is
local so  here $\om$ denotes an open set in $\cd$, and $T$  a (1,1)
positive current in $\om$.

\begin{defi}[\cite{bls}]\label{def_l}  ~
\begin{itemize}
\item[-] $T$ is uniformly laminar if for every $x \in \supp(T)$ there exists 
 open sets $V\supset U\ni x$, with $V$
biholomorphic to the unit bidisk 
$\mathbb{D}^2$ so that in this coordinate chart
$T \rest{U}$ is the direct integral of integration currents over 
a measured family of disjoint graphs in $\mathbb{D}^2$, i.e. :

there exists
a measure $\lambda$ on $\set{0}\times\dd$,  and a family $(f_a)$ 
of holomorphic functions  $f_a: \mathbb{D} \cv \mathbb{D}$ such
that $f_a(0)=a$, 
the graphs $\Gamma_{f_a}$ of two different $f_a$'s are disjoint,  and
\begin{equation}\label{eq_ul}
T \rest{U} = \int_{\set{0}\times \dd}
[\Gamma_{f_a}\cap U ]~ d\lambda (a).
\end{equation}
\item[-] $T$ is laminar in $\om$ if there exists a sequence of open subsets
$\om^i\subset\om$, such that   
$\norm{T}(\fr\om^i)=0$, together with an increasing sequence of
currents $(T^i)_{i\geq 0}$, $T^i$ uniformly laminar in $\om^i$,
converging to $T$. 
\end{itemize}
\end{defi}

Alternatively, a uniformly laminar current is the foliated cycle
associated with an embedded Riemann surface lamination with 
an invariant transverse measure. More generally, we say that 
two disks are  {\em compatible} if they have 
no isolated  intersection points. 
A laminar current always has a {\em laminar
  representation}
\begin{equation}\label{eq_rep}
T= \int_\mathcal{A}[\Delta_\alpha] d\mu(\alpha)
\end{equation}
as an integral over a family of compatible  holomorphic disks, {\em but
with no lamination structure in general}. This means that it is not
possible in general to find open subsets $U$
such that the components of $\Delta_\alpha\cap U$ are closed in $U$.

In a dynamical context, D. Ruelle and D. Sullivan \cite{rs} constructed
uniformly laminar currents subordinate to the stable and unstable
laminations of  a uniformly hyperbolic isolated compact set. Laminar
currents were introduced in \cite{bls} for the purpose of extending
the Ruelle-Sullivan construction to the non-uniformly hyperbolic
setting. The phenomenon of ``folding'', which is apparent in the well
known pictures of the H{\'e}non attractor, is a manifestation of the non uniformity
of the size of the disks in (\ref{eq_rep}).

We proved in \cite{lamin} the following theorem, which gives a very
rough indication of what the local geometry of the Julia sets $J^\pm$
of a general birational map 
should look like. Notice that $\supp(T^\pm)\subset J^\pm$ but whether equality
holds is not known in general. See Diller \cite{di} for definitions
and results related to this question. 

\begin{thm}[\cite{lamin}] If $f$ is an algebraically stable birational map on a
  rational surface $X$ with $\lambda>1$, then the Green currents $T^+$ and
  $T^-$ are laminar.
\end{thm}

The proof is not dynamical: we actually show that any limit of
rational divisors $\unsur{d_n}[C_n]$ in $X$ with 
$$
\text{genus} (C_n) + \sum_{p\in {\rm Sing} (C_n)} n_p (C_n) = O(d_n), 
$$
is a laminar current, and that the Green currents are of this form.
By {\em strongly approximable} we mean a laminar current obtained in this way. 
A crucial point in the present paper 
is that these currents have additional properties,
that were studied in \cite{isect, structure}. We shall discuss many issues
from \cite{isect} in \S \ref{subs_defgeom} and \ref{subs_isect}
below so here we concentrate on \cite{structure}.\\

In the sequel, we will let the dynamics act on the laminar
structure, so we need to know how it is organized. Notice first that the
usual ordering on positive closed currents is compatible with the
laminar structure, in the following sense: if $T_1$ and $T_2$ are
laminar currents with $T_1\leq T_2$, then they admit representations 
$$T_i= \int_\mathcal{A} [\Delta_{\alpha,i}] d\mu_i(\alpha)$$ with
$\mu_1\leq \mu_2$. This  allows us to identify disks and pieces of
laminations subordinate to laminar currents. 

 By definition,
  a {\em flow box} is a closed lamination $\el$ embedded in an open set 
 $U\simeq \dd^2$ such that in 
 this coordinate chart $\el$ is biholomorphic to a lamination by graphs over $\dd$
 ($\dd$ denotes the unit disk). 
 These graphs are called {\em plaques}. If $\el$ is a flow box, we
  define the restriction $T\rest{\el}$ in terms  of the
  representation (\ref{eq_rep}) by integrating only over the disks in 
 $\mathcal{A}$ lying inside one leaf if $\el$.

\begin{defi}[\cite{structure}]\label{def_sub}~
\begin{itemize}
\item[-] A holomorphic disk $\Delta$ is subordinate to $T$ if there
  exists a nonzero uniformly laminar current $S$ with $S\leq T$, and 
 $\Delta$ lies inside a leaf of the lamination associated to $S$. 
\item[-] A flow box subordinate to $T$ is a flow box $\el$ such that
  $\supp(T\rest{\el})=\el$. 
\item[-] The regular set $\mathcal{R}(T)$ is the union of flow boxes,
 or equivalently the union of disks subordinate to $T$.
\end{itemize}
\end{defi}

The main result in \cite{structure} asserts that if $T$ is strongly
approximable, the flow boxes  match
correctly and for every flow box $\el$,
the restriction current $T\rest{\el}$ is uniformly
laminar, i.e. $T$ induces an invariant transverse measure
on $\el$. This will play the role of \cite[\S 4]{bls} in our context.
More precisely, by
{\em weak lamination} we mean a countable union of compatible flow boxes,
 where {\em compatible} here means the associated plaques  do not meet at isolated 
 points.
A {\em  transversal} is by definition a compact set in a flow box
  which meets each plaque at most once. 
 One  feature  of this definition is that the notions of leaf, holonomy, 
 and transverse measure make sense in this setting. 
 
 \begin{thm}[\cite{structure}, Theorems 1.1 and 5.7]\label{thm_structure}
 Let $T$ be a diffuse
 strongly approximable current on the rational surface
 $X$. The regular set $\mathcal{R}(T)$
 has the structure of a weak lamination in the preceding sense. Moreover
 $T$ induces a holonomy invariant transverse measure on this weak lamination. 
 
 If $T$ is extremal as a positive closed current, the transverse measure is ergodic,
 i.e. any measurable saturated set has zero or full measure. 
 \end{thm}
 
The ergodicity of the weak
lamination will be used in the paper through the following
reformulation: 
for any pair of transversals
$\tau_1$, $\tau_2$ of positive transverse measure, there exists a disk
subordinate to $T$ intersecting both $\tau_1$ and $\tau_2$. 
The theorem in \cite{structure} was
stated for $X=\pd$, nevertheless  we explain how to adapt it to a
general rational surface $X$. 

We first need to prove  the fact that $\mathcal{R}(T)$ being a weak
lamination is invariant under birational conjugacy. It suffices to
analyze the action of a birational map $h:\pd \cv X$ on a flow box $\el$. 
Recall that $h$ is the composition of finitely many point 
 blow-ups and inverses of point blow-ups, so it suffices to understand the action 
 of one single blow-up or blow-down on $\el$.  
 
 Let $U$ be an open set such that $\el$ is embedded in $U$.
 If $\pi$ is the blow-up at some point $p\in\el$, 
 $\pi: \pi^{-1} (U)\setminus \pi^{-1}(p)
 \cv U\setminus \set{p}$ is a biholomorphism. Letting 
$\el'$ denote $\el\setminus L(p)$,
where $L(p)$ is the plaque through $p$, 
it is easy to cover $\pi^{-1}(\el')$ with at most countably 
 many flow boxes. The remaining leaf has zero transverse measure so holonomy
  invariance of the transverse measure is not affected.
  
Assume now $T$ is a diffuse strongly approximable laminar current
in $\pi^{-1}(U)$. If $\pi_*T$ 
has non compatible flow boxes, the only possible point of non compatibility is
$p$. But diffuse flow boxes  cannot meet at a single point, so we get 
 a contradiction.

Invariance by holonomy of the induced transverse measure as well as
the  statement concerning ergodicity
 are adapted  in a similar fashion. \hfill $\square$

\begin{rmk}\normalfont 
If $\el$ is a flow box crossing a component $V$ of 
the critical set $C(f)$, then  the images of the plaques 
of $\el$ meet at the point $f(V)$. 
A geometric consequence is the following fact:

{\it For every disk $\Delta$ subordinate to $T^-$ (resp. $T^+$), $\Delta\cap C(f^{-n})$
(resp. $\Delta\cap C(f^n)$) 
is a compatible intersection, that is, either $\Delta\subset C(f^{-n})$ or 
$\Delta\cap C(f^{-n})=\emptyset$.} 

The proof is a simple consequence of the invariance of the currents, together
with Favre's theorem \cite{fa} that 
for every $p\in I(f^n)$, the Lelong number $\nu(p, T^-)$
vanishes.

On the other hand, it is possible for a disk subordinate to $T^+$ to
intersect a component $V$ of $C(f^{-n})$. This will yield a
``pencil'' of plaques  through $f^n(V)$. This phenomenon is seemingly 
observed on computer pictures of stable and unstable manifolds of
birational maps (see e.g. \cite{bed2}). 
\end{rmk}

We will often need to estimate the transverse measure of a given set of plaques. 
If $T$ is strongly approximable, $\el$ is a flow box, and $\tau$ 
is a holomorphic disk
transverse to $\el$, the induced transverse measure on $\tau$ is given by the 
wedge product
$T\rest{\el}\wedge [\tau]$. It is easily proved \cite[Proposition 5.4]{structure}
that if   $\tau$ is the generic (in the measure theoretic sense)
 member of a smooth family of 
holomorphic transversals to $\el$, then the wedge product $ T\wedge
[\tau]$ is admissible (see below \S \ref{subs_defgeom} for a formal
definition)  and
\begin{equation}\label{eq_transv}
 T\rest{\el}\wedge [\tau] =   (T\wedge [\tau])\rest{\el\cap\tau}.
\end{equation}
Abusing notation, if $\tau$ is any transversal 
to the weak lamination (i.e. a closed set 
transverse to a flow box), we will denote the transverse measure induced 
by $T$ on $\tau$ by $T\wedge \tau$. 

\begin{rmk}\normalfont\label{rmk_k3}
Theorem \ref{thm_structure} is precisely  where we use the
rationality assumption on $X$. The dynamical analysis we perform in
the next sections only rely on its conclusions. The invariant currents
associated to automorphisms of projective $K3$ surfaces satisfy these conclusions
(see the remarks in \cite[\S 3]{structure}). In particular 
the discussion to come  also makes sense in that setting, and provide
a new approach to the results in \cite{ca}.
\end{rmk}

 %%%%%%%%%%%%%%%%%%%%%%%%%%%%%%%%%%%%%%%
 %%%%%%%%%%%%%%%%%%%%%%%%%%%%%%%%%%%%%%%
 
\section{Equidistribution of preimages along the unstable current}
\label{sec_equidist}
 
In this section $f$ is an algebraically stable birational map on the
rational $X$, with $\lambda>1$. 
We denote by $\m(\cdot)$ the mass of a current or measure and weak
convergence of currents or measures
is denoted by $\cv$. We normalize invariant currents so
that their mass is 1.
Recall that a transversal 
is by definition a transversal in a flow box. 

The main result in this section is the following equidistribution result. 
It asserts that generic points
on the unstable current $T^-$ become close under backwards
iteration. This approach is new even for complex H{\'e}non mappings. 

\begin{prop}\label{prop_equidist}
If $\tau_1$ and $\tau_2$ are two transversals for the weak
lamination associated to $T^-$, then 
$$(f^{-n})_*\left(\frac{T^-\wedge \tau_1}{\m(T^-\wedge \tau_1)}\right)-
(f^{-n})_*\left( \frac{T^-\wedge \tau_2}{\m(T^-\wedge \tau_2)}\right) \cv 0.$$
\end{prop}

The proposition will be a consequence of the next lemma, which is 
the analogue of the Lyubich-Briend-Duval lemma \cite{ly, brd} in our context. 
Areas  are computed with
respect to the ambient K{\"a}hler form $\omega$.
 
 \begin{lem} \label{lem_lyub}
 Let $\el=\set{D_t, t\in\tau}$ be a flow box subordinate to $T^-$. For every $\e>0$, 
 there exists a positive constant $C(\e)$ and a transversal
 $\tau(\e)\subset\tau$, such that 
 $\m(T^-\wedge \tau(\e))\geq (1-\e) \m(T^-\wedge \tau)$
 and
 $$\forall n\geq 1,~
\forall t\in \tau(\e), ~\area(f^{-n}(D_t))\leq \frac{C(\e)n^2}{\lambda^n}.$$ 
 \end{lem}
 
 \begin{pf} we first analyze the action of the dynamics on the transverse measure.
 Assume $\el$ is a flow box subordinate to $T^-$, avoiding $C(f)\cup I(f)$, and 
 $\tau$ is a transversal in $\el$. Then  
 $f(\el)$ is a flow box for $T^-$ and $f(\tau)$ a transversal in $f(\el)$,
because $f$ is a biholomorphism near $\el$ and $f_*T^-=\lambda T^-$. 
 
 We claim that $T^-\wedge f(\tau)=\lambda^{-1} f_*(T^-\wedge \tau)$. By
 holonomy invariance it suffices
 to prove the result when  $\tau$ lies on a holomorphic disk $\Delta$
 satisfying (\ref{eq_transv}). 
 Then $T^-\wedge \tau= (T^-\wedge \Delta)\rest{\tau}$ is a genuine wedge
 product and
 \begin{equation}\label{eq_inv}
f_*(T^-\wedge\tau)= (f_* T^-)\wedge f(\tau)=\lambda T^-\wedge f(\tau).
\end{equation}
 In particular, since 
  $f_*$ acting on measures preserves masses, this implies 
  $\m(T^-\wedge f(\tau))=\lambda^{-1}\m(T^-\wedge\tau)$. \\
 
 We will now pull back transverse measures. If $\el$ and $\tau$ are as before, 
 moving $\tau$ if necessary we may assume that  $n$ being fixed,
 $\tau\cap C(f^{-n})$ is a finite set of points 
  and $\tau\cap I(f^{-n})=\emptyset$. So up to a set of zero transverse
  measure, $\tau$ is a disjoint union $\tau=\bigcup\tau_j$, with $\tau_j\cap C(f^{-n})=
  \emptyset$. By the previous formula we get that $\m(T^-\wedge f^{-n}(\tau_j))
  =\lambda^n \m(T^-\wedge \tau_j)$, since $f^{-n}(\tau_j)$ avoids $C(f^n)\cup I(f^n)$. 
  
  On the other hand, if $t_1\neq t_2$ in $\bigcup\tau_j$, the disks
   $D_{t_1}$ and $D_{t_2}$ are disjoint and not contained in 
   $C(f^{-n})$ so $f^{-n}(D_{t_1})$ and 
  $f^{-n}(D_{t_2})$ have at most finitely many intersection points. 
  The total $\norm{T}$-mass   of $f^{-n}(\el)$ is 
  $$
  \sum_j \int \area(f^{-n}(D_{f^n(s)})) d(T^-\wedge f^{-n}(\tau_j))(s)\leq \m(T)=1.$$
  Since the total mass of $\sum_j T^-\wedge f^{-n}(\tau_j)$ is $\lambda^n\m
  (T^-\wedge\tau)$, most disks $f^{-n}(D_{f^n(s)})$ have small area
with respect to the transverse measure $\sum_j\left(T^-\wedge
  f^{-n}(\tau_j) \right)$, more precisely
  $$\sum_j\left(T^-\wedge f^{-n}(\tau_j) \right)
  \left(\set{s ,~ \area(f^{-n} (D_{f^n(s)}))
  \geq \frac{cn^2}{\lambda^n}
}\right) \leq \frac{\lambda^n}{cn^2}.$$ Applying $(f^n)_*$ yields
$$(T^-\wedge\tau)\left(\set{t, ~\area(f^{-n}(D_t))\geq \frac{cn^2}{\lambda^n}
}\right) \leq \frac{1}{cn^2}.$$
We now get the conclusion of the lemma by
considering all integers $n$ and adjusting  $c=\frac{\pi^2}
{6\e\m{(T^-\wedge   \tau)}}$. 
\end{pf}
 
 From the lemma we deduce a first equidistribution result. 
 Notice that since transversal
 measures do not charge points, all push forwards $(f^{-n})_*(T^-\wedge \tau)$
 are well defined.
 
 \begin{prop} \label{prop_equidistfb}
 If $\tau_1$ and $\tau_2$ are two global
transversals in a flow box $\el$, then 
 $$(f^{-n})_*\left({T^-\wedge \tau_1}\right)-
(f^{-n})_*\left({T^-\wedge \tau_2}\right) \cv 0.$$
 \end{prop}
 
 \begin{pf} recall from \cite{brd} the following basic area-diameter estimate:
 
\begin{lem} There exists a positive constant $c$, such that if $D\subset 
 \widetilde{D}$ are (possibly singular) disks in $X$, the following
 estimate holds
 $${\rm Diam} (D)^2 \leq c \frac{\area (\widetilde D)}{{\rm Modulus }
 (\widetilde D \setminus D)}.$$
 \end{lem}
  The estimate is only stated for smooth disks in $\mathbb P^k(\cc)$ in \cite{brd},
  however the proof depends only on the Lelong theorem, and the notion 
  of extremal length,  and it carries over for singular disks
  without modification. \\
  
  If $\tau_1$ and $\tau_2$ are closed global transversals in $\el$, they have the 
  same transverse mass by holonomy invariance.  
  Fix a continuous function $\varphi$ on $X$. We must prove 
  \begin{equation}\label{eq_equidist}
  \int \varphi \left[(f^{-n})_*\left({T^-\wedge \tau_1}\right)-
 (f^{-n})_*\left({T^-\wedge \tau_2}\right)\right]\cv 0.
\end{equation} 
 First,  by compactness, there exists a constant $m>0$ such that for every 
 plaque $\widetilde{D}$ of $\el$, there exists a disk $D$, with
  $(\tau_1\cap\widetilde{D})\subset D$ and
 $(\tau_2\cap\widetilde{D})\subset D$ and ${\rm Modulus }
 (\widetilde{D}\setminus D)\geq m$. By the preceding lemmas, 
  for most plaques $\widetilde{D}$, points in $f^{-n}(D)$ get 
  exponentially close under backwards iteration.
   Indeed  for every $\e>0$,  there exists $\tau_i(\e)$, $i=1,2$,  with 
  transverse mass $\m(T^-\wedge \tau_i(\e))\geq (1-\e)\m(T^-\wedge \tau_i)$, such 
  that for $t\in\tau_i(\e)$, ${\rm Diam}(f^{-n}(D_t))^2\leq \frac{cn^2}{m\lambda^n}$.  
  Actually $\tau_1(\e)$ and $\tau_2(\e)$ correspond by holonomy since the property
 that ${\rm Area}(f^{-n}(D_t))$ being small is independent of the transversal.
 
Thus  the term in (\ref{eq_equidist}) writes as 
$$\int \varphi \left[ (f^{-n})_*(T^-\wedge \tau_1\rest{\tau_1(\e)}) -
 (f^{-n})_*(T^-\wedge \tau_2\rest{\tau_2(\e)} )\right] $$
 plus a remainder term not greater than $\e\norm{\varphi}\m(T^-\wedge \tau)$
  because the mass 
 $\m(T^-\wedge \tau_i\rest{\tau_i\setminus\tau_i(\e)})$ is $\m(T^-\wedge \tau)$
  and $(f^{-n})_*$ preserves the mass of measures. 
The latter integral equals
 $$\int_{\tau_1(\e)} \left[ \varphi (f^{-n} (D_t\cap\tau_1)) - 
 \varphi (f^{-n} (D_t\cap\tau_2))  \right] d(T^-\wedge \tau_1)(t)$$ which is small because 
 $\varphi$ is continuous and 
 ${\rm dist} (f^{-n} (D_t\cap\tau_1), f^{-n} (D_t\cap\tau_2))^2\leq \frac{cn^2}{m\lambda^n}$.
 \end{pf}
 
 \noindent{\bf Proof of proposition \ref{prop_equidist}:}
 assume first $T^-\wedge \tau_1$ and $T^-\wedge \tau_2$ have the same 
 (positive) mass. Since $T^-$ is extremal, almost every leaf through $\tau_1$ 
 intersects $\tau_2$ (theorem \ref{thm_structure}). 
 This means that for  $T^-\wedge \tau_1$-almost every point $p$, there 
 exists a disk through $p$, subordinate to $T^-$ and intersecting $\tau_2$. Such a 
 disk is a neighborhood of a path joining $\tau_1$ and $\tau_2$ in the leaf through $p$.
 Fattening those disks in the weak lamination, it
 is standard 
 to prove that for every $\e>0$ there exist  finitely many disjoint
 ``long flow boxes'' $\el_j$, such that $\tau_1\cap\el_j$ and $\tau_2\cap \el_j$ 
 are global transversals in   $\el_j$, and the transverse mass of 
 $\bigcup_j \tau_1\cap\el_j$ and $\bigcup_j \tau_2\cap\el_j$ is greater
 than $(1-\e)\m(T^-\wedge \tau_1)=(1-\e)\m(T^-\wedge \tau_2)$.  Now, 
 the $(f^{-n})_*$ equidistribution of $T^-\wedge \tau_1$ and $T^-\wedge \tau_2$
 follows as in the previous proposition.\\
 
 In the general case choose a large integer $N$. For $i=1,2$, subdivide 
 $\tau_i$ into $E(N\m(T^-\wedge \tau_i))$ pieces $(\tau_{i,j})_j$
 of transverse mass $1/N$, plus  a remainder piece of mass $<\unsur{N}$, 
where $E(\cdot)$ denotes the integer part function. We   
 may moreover assume 
  the measure of $\overline{\tau_{i,j}}\setminus\tau_{i,j}$ is zero.
By the first part of the proof, all pieces $T^-\wedge \tau_{i,j}$ are $(f^{-n})_*$ 
equidistributed, i.e. for any two pairs $(i,j)$ and $(i',j')$,
$$(f^{-n})_*\left({T^-\wedge \tau_{i,j}}\right)-
(f^{-n})_*\left({T^-\wedge \tau_{i',j'}}\right) \cv 0.$$ Thus 
for a continuous function $\varphi$ and  every $i,j$,
$$\limsup_{n\cv\infty} \abs{ 
 \int \varphi (f^{-n})_*({T^-\wedge \tau_i})- E(N\m(T^-\wedge \tau_i))
\int \varphi (f^{-n})_*({T^-\wedge \tau_{i,j}}) } \leq \frac{\norm{\varphi}}{N}. $$ 
This implies, for some constant $c$ depending only on $\m(T^-\wedge\tau_i)$, that
 $$\limsup_{n\cv\infty} \abs{ \unsur{\m(T^-\wedge \tau_1)}
 \int \varphi (f^{-n})_*({T^-\wedge \tau_1})-  \unsur{\m(T^-\wedge \tau_2)}
\int \varphi (f^{-n})_*({T^-\wedge \tau_2}) } \leq \frac{c\norm{\varphi}}{N}.$$
Since $N$ is arbitrary the result follows. \hfill $\square$\\

 %%%%%%%%%%%%%%%%%%%%%%%%%%%%%%%%%%%%%%%
 %%%%%%%%%%%%%%%%%%%%%%%%%%%%%%%%%%%%%%%
\section{The geometric intersection measure}\label{sec_measure}
 
 In order to convert the previous equidistribution statement into a convergence
 result, we need to find an invariant measure with some geometric structure. 
 In this section we  define the geometric intersection 
 of strongly approximable
 laminar currents, and prove that, if non trivial, the geometric intersection measure 
of $T^+$ and $T^-$ has interesting properties. 

\subsection{Geometric intersection}\label{subs_defgeom}
 Geometric intersection of  laminar currents is discussed in 
\cite{bls, isect}.  Nonuniqueness of the laminar representation
(\ref{eq_rep}) makes the general definition of the geometric
intersection measure delicate. In the strongly approximable context,
by using the notion of subordinate disks, we provide a nonambiguous 
definition. \\

Let $T_1=dd^c u_1$ and $T_2= dd^c u_2$ be two closed positive currents
in $\om\subset \cd$. We denote by
$\norm{T}$  the trace measure of the current $T$.
 We say that the wedge product $T_1\wedge T_2$ is
{\em admissible} if $u_1\in L_{loc}^1(\norm{T_2})$. 
Notice that the condition is unambiguous since plurisubharmonic
functions are  defined pointwise. This condition is clearly
independent of the choice of the potential $u_1$ (for convenience
 we drop the {\it loc} subscript).  Under this
condition, the wedge product measure $T_1\wedge T_2$ is defined by 
 $$T_1\wedge T_2= dd^c(u_1T_2).$$ N. Sibony proved (see \cite{these})
that the
admissibility condition is symmetric in $T_1$ and $T_2$ and the wedge
product operation is continuous under decreasing sequences of the
potentials. A useful observation is the following: if
$T_1\wedge T_2$ is admissible and $S_k\leq T_k$, $k=1,2$, are positive
closed currents, then $S_1\wedge S_2$ is admissible and  
$S_1\wedge S_2\leq T_1\wedge T_2$.\\

Following \cite[\S 8]{bls} we now define the geometric wedge product of
uniformly laminar currents.
 
\begin{defi} Let 
$S_1$ and $S_2$ be diffuse uniformly laminar currents, endowed
with representations  $$S_k=\int_{\tau_k}[D_{k,a}]d\mu_k(a),~k=1,2$$
as integrals of families of
submanifolds. We define the product $\geom$ by 
\begin{equation}\label{eq_isectul}
S_1\geom S_2= \int_{\tau_1\times \tau_2}\!
[D_{1,a}\cap D_{2,b}]d\mu_1(a)
d\mu_2(b),
\end{equation}
 where by convention 
the measure $[D_1\cap D_2]$ is the sum 
of Dirac masses at the 
oints of intersection of the disks if they are isolated, zero
if not.
\end{defi}

Since the currents are diffuse, the set of non transverse
intersections has zero measure by \cite[Lemma 6.4]{bls} so counting
multiplicities does not affect the integral in (\ref{eq_isectul}). 
The next proposition asserts that when admissible, $S_1\wedge S_2$
 is
 described as the geometric intersection of the disks constituting 
 $S_1$ and $S_2$ 

\begin{prop}[\cite{isect}, \S 3]\label{prop_isect_ul}
If the wedge product  $S_1\wedge S_2$ is admissible, then $S_1$ and
 $S_2$ have geometric intersection, i.e.
 $S_1\wedge S_2 = S_1\geom S_2$. 

Furthermore, if the leaves of the
 underlying laminations of $S_1$ and $S_2$ only intersect at isolated
 points, then $S_1\wedge S_2$ is admissible.
\end{prop}

We extend the definition of the geometric wedge product
$\dot\wedge$ to sums of 
uniformly laminar currents by summing the geometric intersections of all
factors. We will repeatedly use the obvious fact that the 
 product $\dot\wedge$ is continous under increasing sequences of the factors.

 In the
next proposition, we define a geometric wedge product for 
all strongly approximable  currents.

 \begin{prop}\label{prop_defgeom}
 Let $T_1$ and $T_2$ be two diffuse
 strongly approximable currents on $X$. There exists
 a measure $T_1\geom T_2$, such that if 
 $S_1\leq T_1$ and $S_2\leq T_2$ are uniformly laminar currents in 
$\om'\subset\om$, then
 $S_1\geom S_2\leq T_1\geom T_2$. Furthermore
  $T_1\geom T_2$ has finite mass and
 local product structure (i.e. is a countable sum of product measures).
 
 If the wedge product $T_1\wedge T_2$ is admissible, 
 then $T_1\geom T_2\leq T_1\wedge T_2$.
 \end{prop}
 
 \begin{defi}\label{def_isectgeom}
 If $T_1$ and $T_2$ are two diffuse
 strongly approximable currents on $X$, we say $T_1$ and $T_2$ have 
 non trivial geometric intersection if $\m(T_1\geom T_2)>0$. 
 The measure $T_1\geom T_2$ will be referred to as the geometric intersection
 measure of $T_1$ and $T_2$.
 
 If moreover the wedge product $T_1\wedge T_2$ is admissible and 
 $T_1\geom T_2=  T_1\wedge T_2$, we say that
 $T_1$ and $T_2$ have (full) geometric intersection (or that 
 the wedge product $T_1\wedge T_2$ is geometric).
 \end{defi}
 
 Observe that $T_1$ and $T_2$ have non trivial geometric intersection iff 
 there exist  disks $D_k$, $k=1,2$, respectively subordinate to  $T_k$, with 
 non trivial intersection. 

Recall also that laminar currents were
 defined as increasing limits of currents uniformly laminar in
 $\om_i\subset\om$. Hence if the wedge product $T_1\wedge T_2$ 
is admissible, $T_1$ and $T_2$ have  geometric intersection iff there
 are such increasing sequences $T_k^i \uparrow T_k$, $k=1,2$, 
with $T_1^i\geom T_2^i\cv T_1\wedge T_2$. This is obvious since
 $T_1^i\geom T_2^i \leq T_1\geom T_2\leq
 T_1\wedge T_2$.\\
 
 In order to prove proposition \ref{prop_defgeom} we first give an 
 {\it a priori} bound on 
 masses of geometric intersections. We use hypotheses of global nature. 
 This is actually needed only when the usual wedge product is not
 admissible, 
which is really the new case here.
    
 \begin{lem}
There exists a constant $C$ depending only on $X$
such that if   
 $S_i\leq T_i$,  $i=1,2$, is an at most countable
 sum of uniformly laminar currents $S_i=\sum_j S_{i,j}$, then
 $\m(S_1\geom S_2)\leq C \m(T_1) \m(T_2)$ (in case $X=\pd$ or $\pu\times\pu$, 
the right hand side can be replaced by the intersection pairing
  $(\set{T_1}, \set{T_2})$).

Moreover if the wedge product $T_1\wedge T_2$ is admissible then 
$S_1\geom S_2\leq T_1\wedge T_2$. 
 \end{lem}
 
 \begin{pf} assume first $X=\pd$ or $\pu\times\pu$. In this case we may regularize
 a positive closed current on $X$ by considering a family of
shrinking neighborhoods $(N_\e)$ of $id$
 in ${\rm Aut} (X)$, that is, $T_\e=\unsur{\vol (N_\e)}
 \int_{N_\e} \Phi_* T d\Phi$ is a smooth positive closed
 current, and $T_\e\cv T$ when $\e\cv 0$. Observe that the
 approximation is linear in $T$
 
First, by taking an increasing limit, we may assume that the sums $S_i=\sum_j
 S_{i,j}$ have only finitely many terms (here $i=1$ or $2$). 
 The uniformly laminar current   $S_{i,j}$ is  closed
 in some open set $\Omega_{i,j}$. 
For any small $\alpha>0$, define  $\Omega_{i,j,\alpha}$
 as the open subset 
$$\Omega_{i,j,\alpha}=\set{p\in \Omega_{i,j}, ~ d(p, \fr\Omega_{i,j})\geq \alpha}.$$
As $\alpha\cv 0$, $\Omega_{i,j,\alpha}$ increases to $\Omega_{i,j}$, so
we have the following increasing limit of currents, 
$$S_i= \sum_j S_{i,j}\mathbf{1}_{\om_{i,j}}= \lim_{\alpha\cv 0}
\sum_j S_{i,j}\mathbf{1}_{\om_{i,j, \alpha}}.$$

Fix now $\alpha>0$. For $\e>0$ small enough,
 we can define the regularization
 $S_{i,j,\e}$ in
 $\Omega_{i,j,\alpha}$, and  $S_{i,j}
=\lim_{\e\cv0} S_{i,j,\e}$ in $\Omega_{i,j,\alpha}$.

Define $S_{i, \e}$ as $\sum_j S_{i,j,\e}\mathbf{1}_{\om_{i,j,
    \alpha}}$. We get that 
  $$S_{1,\e}\wedge S_{2,\e}=
  \left(\sum_j S_{1,j, \e} \mathbf{1}_{\Omega_{1,j,\alpha}}\right)\wedge 
  \left(\sum_k S_{2,k,\e} \mathbf{1}_{\Omega_{2,k,\alpha}}
  \right) \leq T_{1,\e}\wedge T_{2,\e}$$ because of the linearity of
  the approximation.
On the other hand the  measure on the right hand side has
  mass $(\set{T_1},\set{T_2})$. 
If the wedge products 
   $S_{1,j}\wedge S_{2,k}$ are locally admissible, for fixed $\alpha$,
    the convergence 
   $S_{1,j, \e}\wedge S_{2,k, \e}\cv S_{1,j}\wedge S_{2,k}$ holds in 
   $\Omega_{1,j,\alpha}\cap\om_{2,k,\alpha}$, when $\e\cv0$. 
 This may be seen for instance as a consequence of 
   geometric intersection of uniformly laminar currents. Hence 
   $$\m\left(S_1\wedge S_2\rest{\bigcup_{j,k}\Omega_{1,j,\alpha}\cap\om_{2,k,\alpha}}
   \right) \leq  (\set{T_1},\set{T_2}). $$
   We conclude that $\m(S_1\wedge S_2)\leq  (\set{T_1},\set{T_2})$  by
   letting  $\alpha$ tend to zero. The second assertion of the 
   lemma is obvious.\\
   
 In the general (non admissible wedge product)
 case just remark that $S_{i,j}$ can be written as  $S_{i,j}=R_{1,j}+
 Q_{i,j}$, where  $R_{1,j}$ is made up of 
 the disks not subordinate to $S_{2,k}$ and having non trivial intersection with 
 $S_{2,k}$. Hence  by definition of the
 geometric wedge product,
 $S_{1,j} \geom S_{2,k}$ equals $R_{1,j}\geom S_{2,k}$. 
 Now the wedge product 
 $R_{1,j}\wedge S_{2,k}$ is admissible by proposition \ref{prop_isect_ul}, so 
 $R_{1,j}\wedge S_{2,k}=S_{1,j}\geom S_{2,k}$ and we conclude as before, 
 replacing $S_{1,j}$ by $R_{1,j}$.\\
   
For an arbitrary rational surface $X$, consider a rational map $h: X\cv
\pd$. Since $T_1$ and $T_2$ are diffuse, $S_1\geom S_2$ 
charge neither points nor curves, so 
$$\m(S_1\geom S_2)=\m\big(S_1\geom S_2
\rest{X\setminus (C(h)\cup I(h))}\big)= \m\big( (h_*S_1\geom h_*S_2)
\rest{\pd\setminus (C(h^{-1})\cup I(h^{-1}))} \big)$$  
where the last equality follows from the fact that 
$h\rest{X\setminus (C(h)\cup I(h))}$ is a biholomorphism. 

Let 
$N$ be any norm on $H^2(X, \cc)$. Observe that 
$(h_*\set{T_1}, h_*\set{T_2})\leq C N(T_1) N(T_2)$, 
because $h_*$ is linear and $H^2(X, \cc)$ is finite dimensional.
 It is
an easy exercice in K{\"a}hler geometry to prove that $N(T_i)$ can be
replaced by the mass $\m(T_i)$ in this inequality.
We now  conclude using the previously discussed case $X=\pd$, because
 $h_*(S_i)\leq h_*T_i$, $i=1,2$, the $h_*T_i$ are strongly approximable
currents on $\pd$. 
\end{pf}
    
 \noindent{\bf Proof of proposition \ref{prop_defgeom}:} 
 Fix a neighborhood basis $(\omega_j)$ of $X$. We assume all $\omega_j$ are 
 biholomorphic to bidisks. For every $\omega_j\simeq\dd^2$, we consider the two 
 sub-bidisks $\omega_j'$ and $\omega_j''$ corresponding to $\bb'=\dd\times D(0,\unsur{4})$
 and $\bb''=D(0, \unsur{4})\times\dd$, where 
 $\dd$ denotes the unit disk. It is a very basic observation that for every line
 $L$
 in $\dd^2$ intersecting $\bb'\cap\bb''$, $L$ is either a graph over the second projection 
 in $\bb'$ or a graph over the first projection in $\bb''$. Rename as $(\Omega_j)$
 the family $(\omega'_j)\cup(\omega''_j)$. 
 
 Let  $D$ be any disk subordinate to $T_1$, and $p\in D$. 
 Because at small scales, $D$ is close to its tangent space at $p$, 
 by the  preceding observation, there exists $\Omega_j\ni p$ such that 
 $D\cap\Omega_j$ is a graph for one of the two natural projections.
 Thus $D=\bigcup_{j \in J_D} D\cap\Omega_j$, where $J_D$ is the set of 
 indices such that $D\cap\Omega_j$ is a graph 
 in the bidisk $\Omega_j$. 
 The open set $\Omega_j$ being  fixed, the set of such graphs in $\Omega_j$ 
 subordinate to $T_1$ form a lamination
 $\el_{1,j}$ in $\Omega_j$. 
 We let $T_{1,j}= T_1\rest{\el_{1,j}}$.
 Doing the same construction with $T_2$, for every $j$ we form the geometric
 intersection measure $\mu_j = T_{1,j}\geom T_{2,j}$.  
 Now the family $\sup (\mu_1, \ldots , \mu_j)$ is increasing and we
 define
 $T_1\geom T_2$ to be its increasing limit, which has finite mass by
  the preceding lemma.\\
 
 Let $S_i\leq T_i$ be uniformly laminar currents.
 Take $p\in\supp(S_1\geom S_2)$, and
 let  $D_1$ be a disk 
 subordinate to $S_1$  through $p$. There exists a bidisk $\Omega^1\ni
 p$ from the neighborhood basis, 
such that $D_1$ is a graph over some direction.
 This also holds for the corresponding leaves close enough to $p$. 
We do the same for $S_2$. Since the 
 disks subordinate to a strongly approximable current are compatible, there is
 at most one disk subordinate to $T_i$ through $p$ so 
 with the preceding notation $S_i\leq T_{i,j}$, and near $p$ 
 in $\Omega^1\cap \Omega^2$,  $S_1\geom S_2\leq T_1\geom T_2$. \\
 
 It only remains to check the product structure. 
 If $\el_i$ is a flow box subordinate to $T_i$, $i=1,2$, then,
 if non trivial, the measure  $T_1\rest{\el_1}\geom 
 T_2\rest{\el_2}$ has product structure. Moreover if $p\in\el_i$, there is exactly one
 disk through $p$ subordinate to $T_i$. Hence 
 $$(T_1\geom T_2)\rest{\el_1\cap\el_2}= T_1\rest{\el_1}\geom T_2\rest{\el_2}.$$ 
 We may now pick a countable collection of 
 disjoint product sets $\el_1\cap\el_2$, 
 of full measure, and $T_1\geom T_2$ has product structure on each of them. 
  \hfill $\square$\\
 
\subsection{Invariant measure} We now study the dynamical properties of the 
geometric intersection measure, provided it is non zero. Recall the
statement of our first main theorem.

\begin{thm}\label{thm_entropy}
Let $f$ be an algebraically stable birational map on a rational
surface $X$, satisfying $\lambda>1$.
 Assume further $T^+$ and $T^-$ have nontrivial geometric
intersection $\mu= T^+\geom T^-$. Then
\begin{enumerate}[\em i.]
\item $\mu$ is an invariant measure which is mixing.
\item  For $\mu$-almost every $p$, there exist unit tangent vectors
$e^u(p)$ and $e^s(p)$ at $p$, depending measurably on $p$, 
there exists $\nn'\subset\nn$ of density $1$,
such that 
\begin{equation}\label{eq_lyap}
\liminf_{\nn'\ni n\cv\infty} \unsur{n} \log\abs{df^n(e^u(p))}\geq 
\frac{\log\lambda}{2} \text{ and } 
\limsup_{\nn'\ni n\cv\infty} \unsur{n} \log\abs{df^n(e^s(p))}\leq 
-\frac{\log\lambda}{2} .
\end{equation}
\item $\mu$ has entropy $h_\mu(f)=\log\lambda$. 
In particular the topological entropy
 $h_{top}(f)$ is $\log\lambda$.
\item $\mu$ has product structure
with respect to local stable and unstable manifolds.
\end{enumerate}
\end{thm}

Item {ii.} requires a few comments. Lyapounov exponents are only defined when
$\log\norm{df}\in L^1(\mu)$. We do not know whether this hypothesis is true
in our context, while (\ref{eq_lyap}) always make sense. 
 Of course when $\log\norm{df}\in L^1(\mu)$ then ii. is a
statement about Lyapounov exponents. 

The bound $\frac{\log\lambda}{2}$ on the Lyapounov exponents is
sharp. 
This is clear if  not only 
birational maps on rational surfaces are allowed, but also holomorphic 
diffeomorphisms on torii: consider for instance the map induced by 
the linear map 
$\left(\begin{smallmatrix}
2&1\\
1&1
\end{smallmatrix}\right)$ on $\cd/\zz[\sqrt{-1}]^2$. Its Lyapounov exponents
--relative  to Lebesgue measure-- are 
$\chi^u= \log\frac{3+\sqrt{5}}{2}>0>\chi^s =
\log\frac{3-\sqrt{5}}{2}$. The topological entropy
is  $2\max(-{\chi^s}, \chi^u)=-2\chi^s$. On the other hand it was
observed by S. Cantat and C. Favre \cite[Example 3.2]{cf} that such a
map gives rise to an automorphism of a {\em rational} surface, 
obtained as the desingularization of  the quotient of the torus
$\cd/\zz[\sqrt{-1}]^2$, by the multiplication by $\sqrt{-1}$. \\

As shown in \cite[p.86]{ow},  {i.} and {iv.} imply $(f,\mu)$ 
has the Bernoulli property, i.e. is measurably
conjugate to a Bernoulli shift.\\ 

\begin{pf} we will prove the items separately. Note that 
  {iv.} follows from proposition \ref{prop_defgeom} as soon as the disks
   subordinate to $T^+$ and $T^-$ are respectively identified as being
   stable and unstable disks, which will be a consequence of the proof
   of {ii.} The measure $\mu$ has finite mass, so by normalization
    we assume  $\mu$ is a probability measure.\\

\noindent{\bf Invariance and mixing.} By hypothesis, $\mu$ is the
   geometric intersection measure of  diffuse laminar currents, so
   $\mu$ gives no mass to subvarieties. In particular we may check the
   invariance of $\mu$ in $X\setminus(I(f^{\pm 1})\cup C(f^{\pm
   1}))$. On any open set $\Omega$ where $f$ is a biholomorphism, 
$f_*(T^+\geom T^-) = f_*T^+\geom f_*T^-=T^+\geom T^-$ so it follows
   that $\mu$ is invariant.\\

The proof that $\mu$ is mixing is slightly reminiscent of the celebrated {\em Hopf
argument} for the ergodicity of the  geodesic flow (see \cite[p. 217]{kh}).  
By construction, $\mu$ is an integral of measures of the form 
$T^-\geom [D]$, where $D$ is a disk subordinate to $T^+$. Moreover
$T^-\geom [D]$ decomposes as
 an at most countable sum of induced transverse
measures $T^-\wedge \tau$ on transversals to $T^-$: indeed this is the
case for every $T^-\rest{\el}\geom [D]$, where $\el$ is a flow box for
$T^-$. 

So by proposition \ref{prop_equidist}, $\mu$ itself is equidistributed with measures of
the form $\frac{T^-\wedge \tau}{\m(T^-\wedge \tau)}$, i.e. 
for every transversal $\tau$, 
$$(f^{-n})_*\mu -   
(f^{-n})_*\left(\frac{T^-\wedge \tau}{\m(T^-\wedge \tau)}\right)
=\mu -   
(f^{-n})_*\left(\frac{T^-\wedge \tau}{\m(T^-\wedge \tau)}\right)\cv 0.$$
If $\varphi$ is a piecewise constant function on a given flow box
$\el$, we get similarly 
\begin{equation}\label{eq_mix}
\mu- (f^{-n})_*\left(\frac{\varphi\mu}{\int \varphi\mu}\right)\cv 0.
\end{equation}
These functions are uniformly dense among
continuous functions on $\el$. Hence (\ref{eq_mix}) holds
 for continuous $\varphi$ on $\el$. 
For global $\varphi$, just write $\varphi=\sum \mathbf{1}_{\el_i}\varphi$, 
where $(\el_i)$ is a collection of disjoint flow boxes of full
$\mu$-measure.
To conclude, we remark that (\ref{eq_mix}) is a reformulation of mixing.\\   

 \noindent{\bf Lyapounov exponents.} 
 We  show that there exists a measurable unit vector field $e^u$, such that
 for fixed $\e>0$, for $\mu$-a.e. $p$, there exists $\nn_\e$ of
 density $\geq 1-\e$ such that 
\begin{equation}\label{eq_lyap2}
\liminf_{\nn_\e\ni n\cv\infty} \unsur{n} \log\abs{df^n(e^u(p))}\geq 
\frac{\log\lambda}{2} .
\end{equation}
It will then suffice to put $\nn'= \bigcup_{\e>0} \nn_\e$. 

Fix $\e>0$ and
 consider a collection $A= \el_1 \cup \ldots \cup \el_N$ 
of disjoint flow boxes for $T^-$, 
 such that $\mu(\el_1\cup\cdots \cup\el_N)\geq 1-\frac{\e}{3}$. 
 If $p\in A$, we denote by $D_p$ the plaque of
  $\el_1\cup\cdots \cup\el_N$ through $p$, and $e^u(p)$ the unit tangent vector
  to $D_p$ at $p$.  
  
  Removing a set of plaques of small $T^-$-transverse measure, hence of small
$\mu$-measure, we get by lemma \ref{lem_lyub} a set still denoted by $A$, with 
measure $\geq 1-\frac{2\e}{3}$, 
such that if $p\in A$, and for every $n$, $\area(f^{-n}(D_p))\leq
\frac{Cn^2}{\lambda^n}$. Making a further reduction we end up with a
set $A$
 with $\mu(A)\geq 1-\e$, such that for each plaque $D$ of 
  $\el_1\cup\cdots \cup\el_N$, $A\cap D$ is relatively compact in  $D$, 
  with a uniform bound on 
  ${\rm dist}(A,\fr D)$.

By the Birkhoff Ergodic Theorem, for a.e. $p$, the 
set 
$$\nn_\e=\set{n\in\nn, ~f^n(p)\in A}$$
has density $\geq 1-\e$. If $p\in A$ and $n\in\nn_\e$ is large enough, 
$f^{-n}D_{f^n(p)}$ has small area, and reducing $D_{f^n(p)}$ slightly
if necessary, small
diameter, so $  f^{-n}D_{f^n(p)}\subset D_p$. Since $D_{f^n(p)}$ lies on
a finite set of flow boxes, by Cauchy estimates (or Koebe distortion), 
the derivative 
$df^{-n}_{f^n(p)}\big(e^u(f^n(p))\big)$ has norm 
$$\abs{  df^{-n}_{f^n(p)}\big(e^u(f^n(p))\big) } = \abs{
\left(df^n_p\right)^{-1}\big(e^u(f^n(p))\big) }\leq 
\frac{Cn}{\lambda^{\frac{n}{2}}}.$$
This gives (\ref{eq_lyap2}) if $p\in A$.

For $\mu$-generic $p$, $p$ does not belong to
 $\bigcup_n C(f^n)$ and for some $n_0(p)$, 
$f^{n_0}(p)\in A$ --more precisely $f^{n_0}(p)$ belongs to the  
full measure subset $A'\subset A$ of points satisfying (\ref{eq_lyap2}). Since 
 $p\notin \bigcup_n C(f^n)$, the differential $df^{n_0}$ is invertible and 
 (\ref{eq_lyap2}) holds at $p$ by pulling back by $f^{n_0}$.\\
 
 We also proved that points in the plaques $D_p$ become exponentially
 close under backwards iteration, so $D_p$ is the local unstable manifold
 of $p$.\\
 
 \noindent{\bf Entropy.} Defining topological entropy requires some care because
 $f$ has indeterminacy points.
 The definition of topological entropy we use is Bowen's definition
 via $(n,\e)$ separated sets on $X\setminus \bigcup I(f^n)$ 
 with respect to  the ambient Riemannian metric. The  Gromov
 inequality \cite{gr}
 asserts that $h_{top}(f)\leq \log\lambda$. The variational principle need not
 hold in this context, 
but the inequality $h_\nu(f)\leq h_{top}(f)$ persists for 
 any invariant probability measure $\nu$. This can be seen for
 instance by restricting to ergodic measures,
 considering partitions by balls of radius $\e$ in the
 definition of metric entropy and applying the Shannon-McMillan-Breiman theorem.\\
 
 Let us recall some material from  entropy theory.
 We shall use the formalism of measurable partitions and conditional measures 
 (see e.g. \cite{bls} for a presentation adapted to our context). 
If $\xi$ is a measurable partition, a probability
 measure $\nu$ may be 
 {\em disintegrated} with respect to $\xi$, giving rise to a probability measure 
 $\nu(\cdot | \xi(x))$ on each atom of $\xi$. We have the following disintegration
 formula: for every continuous function 
$\phi$, 
$$\int \left(\int   \phi(y) d \nu(y|\xi(x)) \right) d\nu(x) = \int\phi d\nu.$$
The partition $\xi$ is said to be $f^{-1}$-invariant if $f^{-1}\xi$ is a refinement
of $\xi$, i.e. for every $x$, $f^{-1}(\xi(f(x)))\subset \xi(x)$. Given partitions
$\xi_i$, we denote by $\bigvee \xi_i$ the {\em joint partition}, i.e.  
 $\left(\bigvee \xi_i\right)(x) 
= \bigcap (\xi_i(x))$. A partition is called a {\em generator} if
 $\bigvee_{n\in\mathbb{Z}} f^{n}\xi$ is the partition into points.
 
 Given a partition $\xi$, we consider the 
 $f^{-1}$-invariant partition $\xi^u=\bigvee_{n\in\mathbb{N}} f^{n}\xi$.
 
 \begin{prop}[Rokhlin] \label{prop_rokhlin}
 If $\xi$ is a generator with finite entropy, then
 $$h_\mu(f) = h_\mu (f, \xi^u)=
  -\int \log \mu(f^{-1}\xi^u(x) | \xi^u(x)) d\mu (x) =
 \int \log J^u_\mu (x)d\mu(x), $$
 where $J^u(x):= \left({\mu\left(  f^{-1}(\xi^u
(f(x)))  |\xi^u(x)\right)}\right)^{-1}$ is the unstable Jacobian. 
\end{prop}

We do not define  the entropy $h_\mu (f, \xi^u)$ here but we stress that the entropy 
 finiteness hypothesis is satisfied because 
 $$ h_\mu (f, \xi^u) \leq h_\mu(f)\leq h_{top}(f)\leq \log \lambda.$$
 
 \begin{prop}[Pesin, see \cite{ls}] \label{prop_ls}
 There exists a measurable $f^{-1}$-invariant generator
 $\xi^u$, whose atoms are open subsets of local unstable manifolds, and 
 such that $h_\mu(f) = h_\mu (f, \xi^u)$.
\end{prop}

\begin{pf} the proposition is stated in the context of Pesin theory in 
\cite[Proposition 3.1]{ls}, 
but it holds in our context. More precisely what is exactly needed in \cite{ls} 
is a family of local unstable manifolds $V_{loc}$ 
satisfying the conclusions
of \cite[Proposition 3.3]{ls}: items (3.3.1) to (3.3.6), except (3.3.5), assert that the 
family of manifolds $V_{loc}$ has controlled geometry on a set of large $\mu$
measure, and (3.3.5) means that points in the same local leaf become exponentially
close under backwards iteration, uniformly on sets of large measure.
The reader will easily check these properties are true for the 
unstable disks constructed above, that is, the set of disks subordinate
to $T^-$. 
\end{pf}

We are now ready to compute $h_\mu(f)$. Consider the unstable partition 
provided by the previous proposition. Since $\mu$
has product structure relative to $T^+$ and $T^-$, for $\mu$-a.e. $x$
$T^+\geom [\xi^u(x)]$ has positive mass.
As an obvious consequence of the product structure of $\mu$
and the definition of geometric wedge product $\geom$, the
conditional measures $\mu(\cdot|\xi^u(x))$ are induced by
$T^+$, more specifically  
$$\mu(\cdot | \xi^u(x)) =
 \frac{T^+\geom[\xi^u(x)]}{\m\left(T^+\geom[\xi^u(x)]\right)}. $$

From the invariance relation $f^* T^+ = \lambda T^+$ (see 
equation  (\ref{eq_inv}))
we deduce that
$$T^+\geom[f^{-1}(\xi^u(f(x)))]= \left(T^+\geom [\xi^u(x)]\right)\rest{ 
f^{-1}(\xi^u(f(x)))} = \unsur{\lambda}f^*(T^+\geom [\xi^u(f(x))]),$$ hence 
the unstable Jacobian $J^u_\mu$ satisfies the multiplicative
cohomological equation 
$$ J^u_\mu (x) = \lambda \frac{\rho(x)}{\rho(f(x))} \text{ a.e. ,  where }
\rho(x) = \m\left(T^+\geom [\xi^u(x)]\right).$$  
Using the invariance of
both $\mu$ and the partition, 
the Birkhoff Ergodic Theorem implies 
$$h_\mu(f) = \int \log  J^u_\mu ~  d\mu = \log\lambda ,$$ see
\cite[Proposition 3.2]{bls} for more details. This concludes the proof of theorem
\ref{thm_entropy}.
\end{pf}

 %%%%%%%%%%%%%%%%%%%%%%%%%%%%%%%%%%%%%%%
 %%%%%%%%%%%%%%%%%%%%%%%%%%%%%%%%%%%%%%%
\section{The Bedford-Diller setting}\label{sec_bd}
 
The aim of this section is to prove that the class of maps considered
in \cite{bd} satisfy the hypotheses of theorem \ref{thm_entropy}.
   The currents $T^+$
and $T^-$ actually have full geometric intersection in this setting. 

\subsection{Geometric intersection}\label{subs_isect}
We prove that under some potential theoretic conditions, the wedge product
of two strongly approximable currents is geometric. The results here
generalize those of \cite{isect} and we refer the reader to this paper
for more details.

The fact that two laminar currents on $X$ intersect
geometrically is a local 
property near every point in $X$. So throughout this paragraph, $\om$ denotes an open 
set in $\cd$. Moreover, reducing $\om$ slightly if necessary, we may replace all 
the $L^p_{ loc}$ conditions by $L^p$ conditions in $\om$.
We first state a local property of strongly approximable currents,
which is proved in \cite[Prop. 4.4]{isect}. 

\begin{prop}\label{prop_mass}
Let $T$ be a strongly approximable laminar current, and $\om\subset\cd$ as above.
 Let $\pi_1$ and $\pi_2$ be generic linear projections.
Then for 
subdivisions $\mathcal{S}_1$, $\mathcal{S}_2$ of the respective projection
bases into squares of size $r$, if 
$$\qq = \set{\pi_1^{-1}(s_1)\cap\pi_2^{-1}(s_2), (s_1, s_2) \in
  \mathcal{S}_1\times \mathcal{S}_2}$$ denotes the associated
subdivision of $\om$ into affine cubes of size $r$, there exists a
current  $T_\qq \leq T$ in $\om$, 
uniformly laminar in each $Q\in \qq$, 
and satisfying the estimate 
\begin{equation}\label{eq_subd}
\m (T-T_\qq)\leq C r^2,
\end{equation}
with $C$ independent of $r$.
\end{prop}

We say that a laminar current satisfying the conclusions of the preceding 
proposition is strongly approximable in $\om$. 

\begin{thm}\label{thm_isect_l2}
Let $T_1=dd^c u_1$ 
and $T_2= dd^c u_2$ be two strongly approximable  
currents in $\om\subset\cd$. Assume $u_1 \in L^1(\norm{T_2})$, 
 $u_2$ has derivatives in $L^2({T_1})$ and 
$u_1$ has derivatives in  $L^2({T_2})$. 
Then the wedge product $T_1\wedge T_2$ is geometric. 
\end{thm}

$L^2$ spaces on positive currents are considered in
\cite{bd,bes}. Let $\om\subset\cd$, and $T$ be a positive
current in $\om$. Let $u,v$ be smooth functions. Following \cite{bd} 
we define the pairing
$$\mathcal{E}(u,v)= \int du\wedge d^cv\wedge T,$$
and denote by $\abs{\cdot}_T$ the associated seminorm, 
$\abs{u}_T=\left(\int du\wedge 
d^cu\wedge T\right)^\unsur{2}$. If $u$ is a p.s.h. function
in $\om$ we say that {\em $u$ has derivatives in} $L^2(T)$ if for every regularizing
sequence $u^j\downarrow u$,  $(u^j)$ is a  Cauchy sequence for 
$\abs{\cdot}_T$. If $u$ has derivatives in $L^2(T)$, then $u$ has 
derivatives in $L^2(S)$ for every $S\leq T$.

\noindent{\bf Proof of theorem \ref{thm_isect_l2}:} we
  follow the approach of \cite{isect} closely, only differing in 
the final estimate. The letter $C$ denotes a constant that may change
from line to line, remaining independent of $r$. 
The currents $T_1$ and $T_2$ being strongly
approximable, by proposition \ref{prop_mass}, there exist
for each $r>0$, a subdivision $\qq$, which we may assume is the same for 
$T_1$ and $T_2$, 
and for each $Q\in \qq$ a uniformly laminar current $T_{k,Q}$, $k=1,2$,
 such that
\begin{equation}\label{eq_subd2}
\m\left( T_k - T_{k,\qq} \right) = 
\m\Bigl(T_k-\sum_{Q\in \qq} T_{k,Q}\Bigl)\leq Cr^2,
~k=1,2.\vspace{-.2em}
\end{equation}
 We have to estimate the mass of 
\begin{equation}\label{eq_mass}
T_1\wedge T_2 \quad  - \sum_{Q \in \qq}
T_{1,Q}\wedge T_{2,Q} ,\vspace{-.2em}
\end{equation} where the second term is a geometric wedge
product  because of uniform laminarity.

The first step is to choose an adapted subdivision so that $T_1\wedge
T_2$ is not too concentrated near the boundary of the cubes. 
More specifically, for
$\lambda<1$ close to 1 and $Q\in\qq$,  let   $Q^\lambda$ be the 
homothetic cube of $Q$ with respect to its center, with factor
$\lambda$. As in \cite[lemma 4.5]{isect}, up to a translation of
$\qq$, we may choose $\lambda$ independent of $r$ so that 
 the mass of $T_1\wedge T_2$ in the union of  
 $Q\backslash Q^\lambda$ is small (i.e. smaller than $2(1-\lambda^4)$).  

To handle the remaining part of (\ref{eq_mass}), $\qq$,  
and $\lambda$ being fixed by now, let $\chi$ be a  
nonnegative $C^\infty$ function, with $\chi=1$ near every
$Q^\lambda$, vanishing near the boundary of
every $Q\in\qq$, and with derivatives bounded by $C/r$ in uniform
norm. The problem reduces  to bounding 
 $$\int \chi \left(T_1\wedge T_2 - T_{1,\qq}\wedge
T_{2,\qq}\right) =
\sum_{Q\in\qq} \int \chi \left((T_1\wedge T_2)\rest{Q} - T_{1,Q}\wedge
T_{2,Q}\right) . $$ Moreover in each cube  $Q$,
$$ T_1\wedge T_2 - T_{1,Q}\wedge T_{2,Q} =
T_1\wedge (T_2 -
T_{2,Q}) + T_{2, Q}\wedge (T_1 - T_{1,Q})
\leq T_1\wedge (T_2 -
T_{2,Q}) + T_2\wedge (T_1 - T_{1,Q}),$$
so  by taking the union over all cubes $Q\in\qq$, we infer that   
$$ T_1\wedge T_2 - T_{1,\qq}\wedge T_{2,\qq}
\leq T_1\wedge (T_2 -
T_{2,\qq}) + T_2\wedge (T_1 - T_{1,\qq}).$$
Of course 
we need only
consider the first term because the hypotheses are symmetric.
Using the Schwarz inequality and (\ref{eq_subd2}), we infer  
\begin{align*}
\int \chi dd^c u_1\wedge (T_2 - T_{2,\qq}) &= - \int d\chi\wedge d^c u_1
\wedge (T_2 - T_{2,\qq})\\
&\leq \left(\int du_1 \wedge d^c u_1\wedge  ( T_2 - T_{2,\qq})\right)^{1/2}  
 \left(\int d\chi\wedge d^c\chi \wedge  ( T_2 -
T_{2,\qq})\right)^{1/2} \\
&\leq \frac{C}{r} \m \left( T_2 - T_{2,\qq}\right)^{1/2}  
\left(\int du_1 \wedge d^c u_1\wedge  ( T_2 - T_{2,\qq})\right)^{1/2} \\
&\leq C \left(\int du_1 \wedge d^c u_1\wedge  ( T_2 - T_{2,\qq})\right)^{1/2} 
=C \abs{u_1}_{ T_2 - T_{2,\qq}}
\end{align*}
where the Stokes theorem is valid because $\chi$ has compact support
and $T_2-T_{2, \qq}$ is closed
in every $Q\in\qq$. 

Let $u^\e_1$ be a regularizing family. We write $u_1= u_1^\e+ (u_1-u_1^\e)$ and use
 the triangle inequality 
 $$\abs{u_1}_{ T_2 - T_{2,\qq}} \leq \abs{u_1^\e}_{ T_2 - T_{2,\qq}}+
 \abs{u_1-u_1^\e}_{ T_2 - T_{2,\qq}} \leq\abs{u_1^\e}_{ T_2 - T_{2,\qq}}+
 \abs{u_1-u_1^\e}_{ T_2}$$
Since $u_1$ has derivatives in $L^2(T_2)$, we may fix $\e$, independent of $r$,
so that 
$\abs{u_1-u_1^\e}_{ T_2}$ is small. For fixed $\e>0$, the function $u_1^\e$
is smooth so 
by weak convergence, 
$\abs{u_1^\e}_{ T_2 - T_{2,\qq}}$ tends to zero 
 when $\m(T_2-T_{2, \qq})$ does, i.e. when $r\cv 0$.
 \hfill $\square$

\begin{rmk}\normalfont\label{rmk_pluripolar} 
Using the same argument together with proposition \ref{prop_ul} allows to prove the following:\\

\noindent{\bf Theorem.} \begin{it} Let $T_1=dd^c u_1$ 
and $T_2= dd^c u_2$ be two strongly approximable  
currents in $\om\subset\cd$. Assume $u_1 \in L^1(\norm{T_2})$, 
 $u_2$ has derivatives in $L^2({T_1})$ and 
$T_1$ gives no mass to pluripolar sets.
Then the wedge product $T_1\wedge T_2$ is geometric. \end{it}\\

\noindent{\bf Proof:} following step by step 
the proof of the previous theorem only allows to prove  that 
$T_2\wedge (T_1-T_{1,\qq})$ tends to zero as $r\cv 0$, i.e. that
$T_2\wedge T_1$ is approximated by the ``semi geometric'' wedge
products $T_2\wedge T_{1,\qq}$. We claim that these wedge products are
geometric. 
 
Indeed, since $T_1$ does not charge pluripolar sets, neither does
$T_{1,\qq}$, and proposition \ref{prop_ul} asserts that in each cube $Q$,
$T_{1,Q}$ is the increasing limit of a sequence 
 of uniformly laminar currents $S_j$ with continuous potential.
Moreover, the 
potentials of  $S_j$ may be chosen to form a  decreasing sequence: just write 
$S_{j+1}=S_j+R_j$ and choose a nonpositive potential for $R_j$. We thus infer that 
 $T_2\wedge S_j \cv T_2\wedge T_{1,Q}$.
 
On the other hand 
$T_2\wedge S_j$ is a geometric wedge product because of theorem \ref{thm_isect_l2}:
indeed a continuous
plurisubharmonic function has derivatives in $L^2(T)$ for any positive
closed current $T$ : this  is a corollary of the 
polarization identity 
$$2 du\wedge d^c u\wedge T = dd^c(u^2)\wedge T - 2udd^cu\wedge T.$$
The theorem is proved.\hfill $\square$\\

Another consequence of proposition \ref{prop_ul} is that if $T_1$ and $T_2$ do
not charge pluripolar sets, then neither does $T_1\geom T_2$. This is the case under
assumption (\ref{eq_bd}) below.
\end{rmk}

\subsection{Dynamics}\label{subs_bd}
 We turn back to the dynamical context, and give the proof of
 theorem \ref{theo_bd}.
 Due to the possibly complicated dynamics of indeterminacy points, it
is not known whether the wedge product $T^+\wedge T^-$ is admissible
in general. E. Bedford and J. Diller \cite{bd} managed to construct
the wedge product measure $\mu= T^+\wedge T^-$ and study some of 
its dynamical properties  under the condition
\begin{equation}\label{eq_bd}
\sum_{n\geq 0} \unsur{\lambda^n} \abs{\log {\rm dist}(f^n(I(f^{-1})),
I(f))}<\infty   
\end{equation} (where ${\rm dist}$ is the ambient
 Riemannian distance function)
which is satisfied for many birational maps, and is symmetric
with respect to $f$ and $f^{-1}$ \cite[Theorem 5.2]{di}.  Under 
this hypothesis, they proved
the following: if $\omega^{+/-}$ are smooth forms representing the cohomology
classes $\theta^{+/-}$, then $T^{+/-}= \omega^{+/-} + dd^c g^{+/-}$,
where $g^{+}$ is a  quasi-p.s.h. function with derivatives in
$L^2(\omega + T^-)$, and similarly  $g^{-}$ has derivatives in 
$L^2(\omega + T^+)$.

In this case, by \cite[\S 3]{bd}, 
 $\mu = T^+\wedge T^-$ is a well defined wedge product,
and $\mu$ has positive mass for cohomological reasons. 
As said before, the wedge product $T^+\wedge T^-$ being geometric is a
 local property near every point in $X$, so 
 by theorem \ref{thm_isect_l2}, 
$\mu=T^+\geom T^-$. 
Hence
theorem \ref{thm_entropy} applies to give the dynamical 
properties of $\mu$--some of which 
(mixing and non zero exponents) were
already given in \cite{bd}. 

As an example, 
Diller shows in \cite[\S 7]{di} that a polynomial birational map in $\cd$,
which is algebraically stable in $\pd$ satisfies condition (\ref{eq_bd}).

Another result in \cite{bd} is that 
 $C^\infty$ functions with logarithmic poles
at points of $I(f)$ are $\mu$-integrable. This is the case in
particular for $-\log {\rm dist}(x,I)$, near $I\in I(f)$, as well as
$\log^+\norm{df}$ and  $\log^+\norm{d^2 f}$. This allows us to use
  the construction of {\em Lyapounov charts} and Pesin's theory 
(see e.g. the appendix of \cite{kh}). 
 
A consequence is the equidistribution of
 saddle orbits, following \cite{bls2}.

\begin{thm}\label{thm_saddle}
Assume that $f$ and $X$ are as in theorem \ref{thm_entropy}, 
and that $\log^+\norm{df}$ and  $\log^+\norm{d^2 f}$ are $\mu$-integrable.

 Then  saddle points are equidistributed towards $\mu$, that is, if
$HPER_n$ denotes the set of saddle periodic points of period $n$, 
\begin{equation}\label{eq_saddle}
\unsur{\lambda^n}\sum_{p\in HPER_n} \delta_p \cv \mu.
\end{equation} Moreover for every $n$
there exists a set $\mathcal{P}_n$ of saddle points with
$\#\mathcal{P}_n/ \lambda^n \cv 1$, such that every $p\in
\mathcal{P}_n$ lies in the support of $\mu$.

Similarly, Lyapounov exponents can be evaluated by averaging
on saddle orbits, that is,
\begin{equation}\label{eq_saddlelyap}
\unsur{\lambda^n}\sum_{p\in HPER_n} \chi^u(p)\cv \chi^u(\mu),
\end{equation}
where $\chi^u(p)$ (resp. $\chi^u(\mu)$) denotes the positive Lyapounov exponent
of $p$ (resp. $\mu$).
\end{thm}

We use the formalism of Pesin boxes from \cite{bls}. Pesin boxes
  are sets $Q$ of positive $\mu$-measure, together with neighborhoods $N(Q)$  so
that for $x\in Q$, $W^s_{loc}(x)$ and $W^u_{loc}(x)$ are transverse
connected boundaryless submanifolds in $N(Q)$.
 Moreover $Q$ may be chosen so 
that the angle between intersecting stable and unstable manifolds in $N(Q)$
is uniformly bounded from below and the resulting $Q$ has product
  structure. In our setting, the measure $\mu$  has product
  structure in Pesin boxes. We denote by $\el^s(Q)$ and $\el^u(Q)$ the
  stable and unstable laminations in $N(Q)$.\\

\begin{pf} following \cite{bls2}, the equidistribution statements
(\ref{eq_saddle}) and (\ref{eq_saddlelyap})
are formal consequences of mixing, product structure and the upper bound
$\lambda^n+C$ on the number of periodic points of period $n$
\cite[Theorem 0.6]{df}. 

It remains to prove that the saddle points constructed with the method
of  \cite{bls2} lie in $\supp(\mu)$; 
let $\mathcal{P}_n$ be this set of saddle points. We adapt the argument of
\cite[\S 9]{bls}. 

Points in $\mathcal{P}_n$ arise as intersection points of stable-like
and unstable-like disks in open neighborhoods $N(Q)$ of Pesin boxes,
biholomorphic to bidisks.   
The important fact is that for any $p\in
\mathcal{P}_n$, there exists a Pesin box $Q$ so that $W^s_{N(Q)}(p)$ is
a global transversal of $\el^u(Q)$ in $N(Q)$, and similarly for   
$W^u_{N(Q)}(p)$. Here the subscript $N(Q)$ means: connected
component of $p$ in $N(Q)$ of the manifold under consideration.
Without loss of generality, we may assume $p$ is a fixed point. 
Consider the restriction currents $T^+\rest{\el^s(Q)}$ and
$T^-\rest{\el^u(Q)}$. These currents have positive mass because
$\mu(Q)>0$ and the leaves of $\el^s(Q)$ (resp. $\el^u(Q)$) are
subordinate to $T^+$ (resp. $T^-$).
Furthermore $\mu\rest{Q} =
T^+\rest{\el^s(Q)}\geom T^-\rest{\el^u(Q)}$. By the hyperbolic Lambda lemma 
(the Inclination lemma), for every leaf $\el^s$ of the stable lamination
$\el^s(Q)$, the sequence of cut-off iterates
$(f^{-n}(\el^s))\rest{N(Q)}$ 
converges in the $C^1$ topology to $W^s_{N(Q)}(p)$. Hence
$$\unsur{\lambda^{n}} (f_{N(Q)}^n)^*
\left(T^+\rest{\el^s(Q)}\right)\leq T^+$$ is a uniformly laminar
current, with leaves arbitrarily $C^1$-close to  $W^s_{N(Q)} (p)$,
where the notation $f_{N(Q)}$ means  all
iterates are successively restricted to $N(Q)$ --this is the Graph Transform
operator for currents.  There is
an analogous result in the unstable direction. In particular if we let
$$ \mu_n = \unsur{\lambda^{2n}}(f_{N(Q)}^n)^*
\left(T^+\rest{\el^s(Q)}\right)  \geom (f_{N(Q)}^n)_*
\left(T^-\rest{\el^u(Q)}\right)$$ then $0<\mu_n\leq \mu$ and the measure 
$\mu_n$ has support arbitrarily close to $p$. 
\end{pf}

\appendix

\section{An alternate approach to theorem \ref{thm_entropy}}
\label{sec_alternate}

The discussion in section \ref{sec_equidist} and \ref{sec_measure}
was designed to avoid the use of Pesin's theory. We sketch here how to
recover theorem \ref{thm_entropy} by allowing Pesin's theory. The point is
to relate the laminar structure of the currents and the stable and
unstable manifolds. This provides yet another
approach to the results in \cite{bls}, \S 4 and 8.

The setting is the following: we adopt the hypotheses of \S
\ref{subs_bd}, that is $f$ is a birational map on $X$ satisfying 
(\ref{eq_bd}). We assume $Q$ is a Pesin box and 
$0<\mu_1 = S^+\wedge S^-\leq \mu$ is a measure supported by $Q$, 
where $S^{+/-}\leq T^{+/-}$ are
uniformly laminar currents in $N(Q)$, and the leaves of the underlying
laminations are disks.  For $x\in Q$,  we let $S^{+/-}(x)$ be the
disk of the corresponding current $S^{+/-}$ through $x$. 

\begin{prop}\label{prop_appendix} With notations as above,
for $\mu_1$ a.e. $x\in Q$, $S^+(x)\subset W^s_{loc}(x)$ (resp. 
$S^-(x)\subset W^u_{loc}(x)$).
\end{prop}

Before proving the result, we make two observations.  
The first is that
a current $T=dd^cu$ with $du\in L^2_{loc}$
  gives no mass to pluripolar sets. This is classical for measures  
in $\cc$ \cite[Theorem III.7]{t}, and easily extends to currents by
slicing.
The second observation is the following useful proposition 
\cite[Prop. 6.1]{structure}.

\begin{prop}\label{prop_ul}
Let $S$ be a uniformly laminar current, integral of holomorphic graphs
in the bidisk, $S = \int [\Gamma_\alpha] d\mu(\alpha)$. Assume $S$
 gives no mass to pluripolar sets. Then  $S$ can be written as a countable
sum $S=\sum
S_j$, where the $S_j =  \int [\Gamma_\alpha] d\mu_j(\alpha)$ have
continuous potential and disjoint support. 
\end{prop}

\noindent{\bf Proof of proposition \ref{prop_appendix}:} 
we prove the result for $S^+$. One may assume from the
previous observations that $S^+$ and $S^-$
have continuous potentials.

Suppose the result is false,
that is, there exists $R\subset Q $ of positive $\mu_1$-mass such that for
$x\in R$, $W^s_{loc}(x)\neq S^+(x)$. Slightly moving $x$ if necessary,
makes the intersection between $W^s_{loc}(x)$ and $S^+(x)$
 transverse. Indeed, $\mu_1$ has product structure with
respect to $S^+$ and $S^-$, and $R$ has positive measure, so we may
assume the $([S^+(x)]\wedge S^-)$-mass of $R$ inside $S^+(x)$ is
positive. For $y\in R\cap S^+(x)$ near $x$, $W^s_{loc}(y)$ is
transverse to $S^+(y)=S^+(x)$ since the local stable manifolds are disjoint
(see \cite[Lemma 6.4]{bls}). Without loss of generality we write $x$
for $y$. Reducing $N(Q)$ once again, we assume $S^+(x)$ is a global
transversal to the family of stable manifolds. This does not affect
the fact that $\mu_1(R)>0$.  

A corollary of transversality is that $S^+\wedge [W^s_{loc}(x)]>0$
because a set of positive transverse measure of disks intersect
$W^s_{loc}(x)$ transversally --the existence of the wedge product is
ensured since $S^+$ has continuous potential. \\

On the other hand, the current $S^-$ induces the measure $S^-\wedge
[S^+(x)]$ on the disk $S^+(x)$, which is a measure with continuous
potential on $S^+(x)$. Up to a normalizing factor this
measure coincides with the conditional 
measure $\mu_1(\cdot | S^+(x))$. Now $$\mu_1 (\cdot
| S^+(x))\rest{R\cap S^+(x)}\leq \mu_1(\cdot | S^+(x))$$  so the
restriction 
$\nu= \mu_1 (\cdot
| S^+(x))\rest{R}$ has continuous potential also. 

Let $C=\int [W^s_{loc}(y)] d\nu(y)$ be the uniformly laminar current
constructed from the transversal $S^+(x)$ to the family of stable
manifolds and the conditional 
measure $\nu$. One proves (see \cite[Lemma 6.4]{structure}) that $C$
has continuous potential. From the previous discussion and
geometric intersection it follows that $C\wedge S^+>0$. 

Let $\psi$ be a nonnegative test function in $N(Q)$, $\psi=1$ near
$x$. Then $$0<\int \psi C\wedge S^+\leq \int \psi C\wedge T^+.$$
Consider now
the sequence of currents $\lambda^{-n} (f^n)_*(\psi C)$, where the action
$f_*$ has to be understood here as a proper transform near indeterminacy 
points. A result of J. Diller \cite{di}, built on classical arguments,
asserts that the cluster points of this sequence of currents are
positive closed currents of mass $\int \psi C\wedge T^+>0$. However,
$C$ is an integral of local stable manifolds so $\m (\lambda^{-n}
(f^n)_*(\psi C) )\cv 0$. We have reached a contradiction. \hfill $\square$\\

\bigskip

\noindent{\sc \small UFR de math{\'e}matiques,
Universit{\'e} Paris 7,
Case 7012,
2 place Jussieu,
75251 Paris cedex 05, France.}\\
{\tt dujardin@math.jussieu.fr}

\end{document}